\documentstyle[12pt]{article}

\def\dsp{\displaystyle}

\def\r#1{(\ref{#1})}

\def\bel{\begin{equation}\label }
\def\ee{\end{equation}}
\def\ba{\begin{array}}
\def\ea{\end{array}}

\def\grad{\mathop{\rm grad}}

\begin{document}
\def\tens{\otimes}
\def\ben{\begin{enumerate}}
\def\een{\end{enumerate}}
\def\g#1#2{\grad(#1,#2)}
\def\hg{\hat g}
\def\<{\left<}
\def\>{\right>}
\title{NEW QUANTUM DEFORMATION OF $OSP(1;2)$
\thanks{Talk given at the 7-th Colloquium "Quantum Groups
and Integrable Systems", Prague 1998.}}
\author{C. JUSZCZAK and J. T. SOBCZYK
}
\date{\normalsize
Institute of Theoretical Physics, University of Wroc\l{}aw,\\
\normalsize Pl. Maksa Borna 9,
50-204 Wroc\l{}aw, Poland\vspace*{-3mm}}
\maketitle
\begin{abstract}
The complete set of formulas describing the new quantum deformation of
the $OSP(1;2)$ supergroup is provided. A general ansatz is solved for 
the
deformation of the Borel subalgebra of its dual quantum deformation of
$osp(1;2)$.
\end{abstract}
\newpage
\section{Introduction}\label{sec1}

In the recent paper \cite{Juszsob} all the Lie
superbialgebra structures for the Lie
superalgebra $osp(1;2)$ were classified. All of them turned out to be
coboundaries.
The complete list of inequivalent classical $r$ matrices consists of 
three
entries, one being a one-parameter family.
A natural question is what are
corresponding quantum deformations of both the universal
enveloping superalgebra $U(osp(1;2))$ and of the supergroup 
$OSp(1;2)$.
In the cases of Lie superbialgebra structures $r_1$ and $r_3$ 
everything
is known
\cite{Kulish}, \cite{KCh},
\cite{CK}. $r_2$ has been studied in the recent paper of
Kulish \cite{Kultwist}. The idea was to find a twist element giving 
rise
to triangular Hopf algebra. Unfortunately even if some information 
about
the twist was deduced, its final form remains yet unknown. Therefore 
the
problem of finding the deformation of $U(osp(1;2))$ in the direction 
of
$r_2$ remains open. What one can do is to provide a complete set
of formulas defining the dual quantum supergroup $OSP_p(1;2)$.
Because of relevance of $OSP(1;2)$ supergroup \cite{Sal} it is of
potential
physical interest.

The paper is organized as follows. In Sec. II basic definitions and
notations are introduced (\cite{Juszsob} is followed closely).
Sec. III contains very few details about
the computation techniques used as they
are standard \cite{FRT}. Then
all the relations defining the quantum supergroup
$OSP_p(1;2)$ are presented.
In Sec. IV some problems related with obtaining formulas defining
quantum deformation
of $U_p(osp(1;2))$ are discussed.

\section{Basic definitions and notation}\label{sec2}

The Lie superalgebra $osp(1,2)$ is spanned by the set of generators
$\{ H$, $X_+$, $X_-$, $V_+$, $V_-  \}$
of grading
$grad (H) = grad (X_{\pm})=0$ and $grad (V_{\pm})=1$
which fulfill the following (anti)-commutation
relations
\bel{8a}
\ba{rclrclrcl}
[H, X_\pm] &=& \pm X_\pm\,,\qquad &[H, V_\pm] &=&
\pm \frac12 V_\pm\,,\\{}
[X_+, X_-] &=& 2H\,,\qquad & [X_\pm, V_\pm] &=& 0\,,\qquad \\{} 
 \{V_+,V_-\} &=& -\frac12 H\,,& [X_\pm, V_\mp] &=& V_\pm\,, \\{}
\{V_\pm,V_\pm\} &=& \pm\frac12 X_\pm\,.\qquad 
\ea
\ee
For the supergroup $OSp(1|2)$ we will use the
parameterization by means of $3\times 3$ supermatrices 
\cite{Kulish}

\bel{ag}
T = \pmatrix{a&\alpha&b\cr \gamma&e&\beta\cr c&\delta&d}
\ee
subject to constraints:
$e=1+\alpha\delta$, \ $\gamma = c\alpha - a\delta$, \
$\beta = d\alpha - b\delta$, \  $ad - bc + \alpha\delta =1$.
The variables denoted by Greek letters are of Grassmanian type.
Coproducts follow from matrix multiplication of the elements of $G$.

While the general expression for the antisymmetric $r$-matrix takes 
one of the two forms:
\bel{15}
r_{I} =  x^2 H\wedge X_+ + xy X_+\wedge X_- + y^2 H\wedge X_- \,.
\ee
\bel{15a}
r_{II} = x (H \wedge X_+ - V_+ \wedge V_+) +
         y (X_+ \wedge X_- + 2 V_+ \wedge V_-)
+ z (H \wedge X_- - V_- \wedge V_-)
\ee
($x$, $y$, $z$ -- being arbitrary complex numbers), 
only three nonequivalent $r$-matrices can be singled out:
\bel{15b}
r_1\equiv H\wedge X_+\,,
\ee
\bel{14}
r_2\equiv H \wedge X_+ - V_+ \wedge V_+\,,
\ee
\bel{14a}
r_3 \equiv   t(H \wedge X_+ - V_+ \wedge V_+ +H \wedge X_- - V_- 
\wedge
V_-)\,.
\ee

The parameter $t$ in $r_3$ becomes irrelevant upon quantization as it 
can 
be absorbed into a deformation parameter; it cannot however be 
removed by a change of the basis in the Lie super-bialgebra.

By applying the automorphism (change of basis) $r_3$ can be 
rewritten in the equivalent form as
\bel{14c}
r_3=\tilde t (X_+\wedge X_- + 2V_+\wedge V_-)\,.
\ee

In order to express it in the way in which it is present in
\cite{Kulish} it is necessary to add the symmetric ad-invariant 
combination
\bel{14d}
\tilde t (2H\otimes H + X_+\otimes X_- + X_-\otimes X_+ +
2(V_+\otimes V_- - V_-\otimes V_+)).
\ee
The classical $r$-matrix $r_1$ gives rise the the triangular 
deformation of $osp(1;2)$ discussed in \cite{CK}.

In this paper a correspondence
between (graded) tensor products and (super) matrices is explored.
$3\times 3$ supermatrices are graded in the following way

\bel{j1}
grad(a_{jk}) = \delta_{j2} +\delta _{k2}.
\ee

The correspondence mentioned above is given by \cite{Pittner}

\bel{j2}
e_{ij}\otimes e_{kl} \longrightarrow (-1)^{(g(i)+g(j))g(k)}
E_{3(i-1)+k,3(j-1)+l}\,.
\ee

With the help of above identification and matrix representation of
generators of $osp(1;2)$:

\bel{mat}
H={1\over 2}\pmatrix{1&0&0\cr0&0&0\cr 0&0&-1},\
X_+=\pmatrix{0&0&1\cr 0&0&0\cr 0&0&0},\
V_+={1\over 2}\pmatrix{0&1&0\cr 0&0&1\cr 0&0&0}
\ee

the $r_2$  can be represented as

\bel{j3}
r_2={1\over 2}\left( \begin{array}{ccc|ccc|ccc}
                       0&0&1&0&-1&0&-1& 0&0\\
                       &0&0&0&0 &1&0 & 0&0\\
                       &&0&0&0 &0&0 & 0&1\\
                       \hline
                       &&&0&0 &0&0 &-1&0\\
                       &&&&0 &0&0 & 0&1\\
                       &&&& &0&0 & 0&0\\
                       \hline
                       &&& &&& 0 & 0 &-1\\
                       &&& &&& & 0 & 0\\
                       &&& &&&  &  & 0
            \end{array}\right)           
\ee

Quantum $R$ matrix satisfying (graded)
YB equation can be calculated as $R=\exp({2p\,r_2})$ (the factor 2 is 
taken 
for convenience). One obtains

\bel{j4}
R = \left( \begin{array}{ccc|ccc|ccc}
               1 & 0 & p & 0 &-p & 0 &-p & 0 & {p^2\over 2} \\
               & 1 & 0 & 0 & 0 & p & 0 & 0 & 0\\
               &   & 1 & 0 & 0 & 0 & 0 & 0 & p\\
               \hline
               &   &   & 1 & 0 & 0 & 0 &-p & 0\\
               &   &   &   & 1 & 0 & 0 & 0 & p\\
               &   &   &   &   & 1 & 0 & 0 & 0\\
               \hline
               &   &   &   &   &   & 1 & 0 &-p\\
               &   &   &   &   &   &   & 1 & 0\\
               &   &   &   &   &   &   &   & 1
             \end{array} \right)  .
\ee

One should mention that the above $R$ can also be obtained by applying
some similarity transformation to $R$ given in \cite{Kulish} and then 
by
performing the suitable limit \cite{Kultwist}.

An $\tilde R$
matrix satisfying (not graded) YB equation can be obtained from \r{j4}
by means of \cite{Basz}

\bel{gyb}
\tilde R_{ijmn} = (-1)^{grade(i)\cdot grade(m)} R_{ijmn}
\ee

\section{Quantum $OSp_p(1;2)$ supergroup}

Basic ingredients of the $RT_1T_2=T_2T_1R$ formalism of \cite{FRT}
are $9\times 9$ matrices:

\bel{j5}
T_1 = T\otimes 1= \left( 
                    \begin{array}{ccc|ccc|ccc} a&&&\alpha&&&b&&\\
                              &a&&&-\alpha&&&b&\\
                              &&a&&&\alpha&&&b\\
                              \hline
                              \gamma&&&e&&&\beta&&\\
                              &-\gamma&&&e&&&-\beta&\\
                              &&\gamma&&&e&&&\beta\\
                              \hline
                              c&&&\delta&&&d&&\\
                              &c&&&-\delta&&&d&\\
                              &&c&&&\delta&&&d
                    \end{array}\right)
\ee

and

\bel{j6}
T_2 = 1\otimes T = \left( \begin{array}{ccc|ccc|ccc}
                             a&\alpha& b&&&&&&\\
                               \gamma&e&\beta&&&&&&\\
                               c&\delta&d&&&&&&\\
                               \hline
                               &&&a&\alpha& b&&&\\
                               &&&\gamma&e&\beta&&&\\
                               &&&c&\delta&d&&&\\
                               \hline
                               &&&&&&a&\alpha& b\\
                               &&&&&&\gamma&e&\beta\\
                               &&&&&&c&\delta&d
                        \end{array} \right)       
\ee

$RT_1T_2=T_2T_1R$ can be rewritten as a set of $81$ quadratic 
equations
for entries $a, \alpha, b, ...$.
They have to be supplemented by deformed superorthogonality condition:
\bel{j7}
TCT^tC^{-1} = CT^tC^{-1}T = 1
\ee
where the supermatrix $C$ satisfies
\bel{j8}
R = C_1(R^{t_1})^{-1}C_1^{-1}.
\ee

A form of $C$ can be identified as

\bel{j9}
C=\pmatrix{p&&-1\cr &1&\cr 1&&}.
\ee

It is now possible to write down the complete set
 of relations defining the new
quantum deformation of $OSp(1;2)$. The deformation parameter is $p$.

\bel{wzorki}
\begin{array}{lll}
 [a, b] = p (a^2-1),& [a,c] =-pc^2, & [a,d] = p(ca -cd),\\{}
& [b,c] = -p ca - pdc, & [b, d] = p (1-d^2),\\{}
 & & [c,d] = pc^2, \\{}
[a, \alpha] =0, & [a,\delta ] = -pc\delta , \\{}
[b, \alpha ] = -p\alpha a, & [b,\delta ] = - p(d\delta + c\alpha 
),\\{}
[c,\alpha ] = c\delta, & [c,\delta ] = 0, \\{}
[d,\alpha ]= p(\delta d - \delta a), &[d, \delta ] = -p\delta c, \\{}
\alpha^2 = {p\over 2} (1-a^2), &
\{\alpha ,\delta \}= p(\delta^2-ac), &\delta^2 = -{p\over 2} c^2. \\{}
\end{array}
\ee

It is sufficient to concentrate only on the
relations containing $a, b, c, d, \alpha ,\beta$ as
remaininig ones are not independent.
In order to demonstrate that one first has to
solve for $e$ the square root equation 
$e^2=1+2\alpha\delta+pac-{p\over 2}\delta^2$
obtaining
\bel{s1}
e = 1 + \alpha\delta + {p\over 2} ac.
\ee
Its inverse turns out to be
\bel{s2}
e^{-1} = (1 - \alpha\delta - {p\over 2}ac){1 \over 1 -{1\over 
4}p^2c^2}.
\ee
From the relation $a\delta + \gamma e - c\alpha - {p\over 2} c\delta 
=0$
one derives
\bel{s3}
\gamma = \alpha c - \delta a + p \delta c\,.
\ee
Similarly one calculates
\bel{s4}
\beta = \alpha d-\delta b +{p\over 2}\alpha c - {p\over 2} \delta a
+p\delta d +{p^2\over 2} \delta c\,.
\ee
Finally, it is straithforward to get the expression for the antipode
\bel{ant}
S\pmatrix{a&\alpha&b\cr \gamma&e&\beta\cr c&\delta & d}
=\pmatrix{d-{p\over 2}c&-\beta +{p\over 2}\gamma&-b+{p\over 2}(a-d) +
{p^2\over 4}c\cr
\delta&e&-\alpha -{p\over 2}\delta\cr
-c&\gamma&a+{p\over 2}c}\,.
\ee

\section{Duality}

In an attempt to obtain the dual quantum deformation of the
universal envelopping
algebra of $U(osp(1;2))$ one faces a problem that both $R^{(+)}$ and
$R^{(-)}$ are upper triangular supermatrices
and apparently one cannot get more than
a deformation of Borel subalgebra of $osp(1;2)$ generated by $H, X_+$
and $V_+$.

If the entries of $L_+$ matrix are denoted as
\bel{dl}
L^+ =\pmatrix{A&B&C\cr 0&1&E\cr 0&0&F}
\ee
then $R^+L^+_1L^+_2=L^+_2L^+_1R^+$ equations are
\bel{wzorki2}
\begin{array}{ll}
[A, C]= p(AF-A^2)\,,\qquad & [C, F]=p(F^2-AF)\,, \\{}
[A, F]=0\,,\hfill [C, A] -\!\! &[C, F] = E^2-{p\over 2}F^2
   +{p\over 2}A^2 +B^2\,,\\[2mm]
[A, B] = 0\,, \qquad & [B, F]=0\,, \\{}
[A, E] = 0\,, \qquad & [E, F]=0\,, \\{}
[B, C] = p(BF -BA-E)\,,\qquad & [C, E]= p(FE+B-AE)\,,\\[2mm]
\multicolumn{2}{l}{
B^2={p\over 2} (A^2-1)\,,\qquad \{B, E\} = p(A-F)\,,\qquad
E^2={p\over 2} (1-F^2)\,.}
\end{array}
\ee

One can try to solve them by making an ansatz
\bel{ansatz}
\begin{array}{lll}
A=K(X_+)\,,\quad& F=L(X_+)\,,\\{}
B=V_+\cdot M(X_+)\,,\quad& E=V_+\cdot N(X_+)\,,\quad &C=H\cdot 
P(X_+)\,,
\end{array}
\ee
where $H$, $X_+$ and $V_+$ satisfy undeformed (anti-)commutation 
relations.
$K$, $L$, $M$, $N$ and $P$ are arbitrary analytic functions in one 
variable.
One obtains a complicated set of equations which  lead to 
the following conditions:
\bel{czarek}
\begin{array}{ll}
K\cdot L=1, & M=K\cdot N\,,\\[2mm]
P= \dsp {p(K^2-1)\over X_+\cdot K'}\,,\qquad& 
\dsp {X_+\over 2}N^2=\dsp p {K^2-1\over K^2}\,.
\end{array}
\ee
It is interesting to observe that expressions obtained in
\cite{Kultwist} provide a particular solution of the above set of
equations, namely
\bel{czarek2}
\begin{array}{rcl}
K(X_+)&=&e^{\sigma}\equiv (1+pX_+)^{1\over 2}\,,\\
P(X_+)&=&2pe^{\sigma}\,,\\
N(X_+)&=&\sqrt{2}pe^{-\sigma}\,,\\
M(X_+)&=&\sqrt{2}p\,.
\end{array}
\ee
In such a way one gets a deformation of the Borel subalgebra of
$osp(1;2)$ with
standard (undeformed) algebraic relations and coproducts
\bel{kopy}
\begin{array}{l}
\Delta (e^{\sigma}) = e^{\sigma}\otimes e^{\sigma}\,,\\{}
\Delta (V_+) = e^{\sigma} \otimes V_+ + V_+\otimes {\bf 1}\,,\\{}
\Delta (H) = {\bf 1}\otimes H + pV_+e^{-\sigma}\otimes
V_+e^{-2\sigma}+ H\otimes e^{-2\sigma}\,.
\end{array}
\ee
The coproduct for $X_+$ can be expressed in equivalent form as
\bel{kopy2}
\Delta (X_+) = X_+\otimes{\bf 1}+ {\bf 1}\otimes X_++ pX_+\otimes
X_+\,.
\ee
The calculation of the complete set of formulas describing
$U_p(osp(1;2))$ turns out to be more complicated \cite{Ohn}.
One idea is to calculate the twist element and its form was partially
deduced in
\cite{Kultwist}. Other idea is to apply a suitable limiting procedure
in the $L^{(\pm)}$ formalism of \cite{FRT}. This problem is now under
consideration.

\subsection*{Acknowledgments}

The authors (partially supported by KBN grant 2P 30208706)
would like to thank P. Kulish for valuable discussions.


\begin{thebibliography}{99}
\frenchspacing
\bibitem{Juszsob}{C. Juszczak, and J. T. Sobczyk, 
{\em Classification of low dimensional Lie super-bialgebras}, 
q-alg/9712015, to be published in J.
Math. Phys. 1998.}

\bibitem{FRT}{L. D. Faddeev, N. Yu. Reshetikhin and L. A. Takhtajan,
Leningrad Math. J. {\bf 1} (1990) 193.}


\bibitem{Sal}H. Saleur, Nucl. Phys. {\bf B336} (1990) 363.

\bibitem{Kulish}P. P. Kulish, {\sl
Quantum superalgebra $osp(2|1)$ }, Preprint RIMS-615,
Kyoto, 1988; P. P. Kulish and N. Yu. Reshetikhin,
Lett. Math. Phys. {\bf 18}
(1989) 143.

\bibitem{KCh}M. Chaichian and P. P. Kulish, {\sl Quantum Group
Covariant

Systems}, in "From Field Theory to Quantum Groups", birthday volume
dedicated to Jerzy Lukierski, eds. B. Jancewicz and J. Sobczyk,
World Scientific, 1996, pp. 99-111.

\bibitem{CK}E. Celeghini and P.P. Kulish, J. Phys. A {\bf 31} (1998) 
L79.

\bibitem{Pittner} L. Pittner, P. Uray, J. Math. Phys. {\bf 36} (1995) 
944.


\bibitem{Ohn} Ch. Ohn, Lett. Math. Phys. {\bf 25} (1992) 85.


\bibitem{Basz} V. V. Bazhanov and A. G. Shadrikov,
  Theor. Math. Phys. {\bf 73} (1987) 402.

\bibitem{Kultwist}  P. P. Kulish, {\em Super-jordanian deformation 
of the orthosymplectic Lie superalgebras}, math.QA/9806104.

\end{thebibliography}
\end{document}